\documentclass[12pt]{amsart}
\usepackage{amsmath}
\usepackage{amsthm}
\usepackage{amssymb}
\usepackage{amscd}

\newcommand{\brq}{^{[q]}}

\newcommand{\deff}[1]{{\sl #1}}

\newcommand{\inc}{\subseteq}
\newcommand{\Ga}{\Gamma}
\newcommand{\gothic}{\mathfrak}
\renewcommand{\hat}{\widehat}

\newcommand{\length}{\lambda}
\newcommand{\limt}{\lim\limits_{\longrightarrow t}}
\newcommand{\m}{{\gothic{m}}}

\renewcommand{\phi}{\varphi}

\newcommand{\Rq}{{R^{1/q}}}

\renewcommand{\to}{{\longrightarrow}}
\newcommand{\Tor}{\operatorname{Tor}}

\newtheorem{thm}{Theorem}

\newtheorem{lemma}[thm]{Lemma}

\begin{document}
\title{The vanishing of $\Tor_1^R(R^+,k)$ implies that $R$ is regular}
\author{Ian M. Aberbach}
\address{Mathematics Department \\
 	University of Missouri\\
	Columbia, MO 65211 USA}
\email{aberbach@math.missouri.edu}
\urladdr{http://www.math.missouri.edu/people/iaberbach.html}

\date{\today}
\thanks{The author was partially supported by  the National Security Agency.
He also wishes to thank the referee for a careful reading of this paper
and several corrections.}

\bibliographystyle{amsplain}


\begin{abstract} Let $(R,m,k)$ be an excellent local ring of positive prime
characteristic.  We show that if $\Tor_1^R(R^+,k) = 0$ then $R$ is 
regular.  This improves a result of Schoutens, in which the additional
hypothesis that $R$ was an isolated singularity was required for the proof.
\end{abstract}  

\maketitle

Let $R$ be an integral domain.  Then we denote by $R^+$ the integral
closure of $R$ in an algebraic closure of the fraction field of $R$.
Under the assumption that $R$ is a local excellent domain with positive prime characteristic $p$, the ring
$R^+$ is a balanced big Cohen-Macaulay algebra \cite{HH:1992}.  We assume for
the rest of this paper that $R$ is a commutative ring with positive prime
characteristic $p$.  Let $F:R \to R$ be the Frobenius
endomorphism given by $r \mapsto r^p$.  It is a theorem of Kunz \cite{Kunz:1969} that $R$ is regular if and only if $F$ is a flat map.  From this theorem it
is not difficult to show that $R$ is regular if and only if $R^+$ is flat
over $R$.  The more general question of whether $\Tor_1^R(R^+,k) =0$ implies
that $R$ is regular for a local ring $(R,\m,k)$ of positive characteristic
is posed in the exercises
in section 8 of \cite{Huneke:1996}  (when $\Tor_1^R(S,k) =0$ for a module-finite
extension then Nakayama's lemma shows that $S$ is flat over $R$, however,
$R^+$ is far from finitely generated over $R$).  Schoutens has shown that
for an excellent local ring the condition $\Tor_1^R(R^+,k) =0$ implies that $R$ is weakly $F$-regular,
and if $R$ has an isolated singularity then $R$ is regular (\cite{Schoutens:2003}, Theorems 1.3 and 1.1).  We show here that, in fact,
the vanishing of $\Tor_1^R(R^+,k)$ suffices to imply regularity for excellent
rings of positive prime characteristic.

Assume that $(R,\m,k)$
is a reduced excellent local ring.  $R$ is then approximately Gorenstein,
so there is a sequence of irreducible $\m$-primary ideals $\{I_t\}$ cofinal
with the powers of $\m$ (see \cite{Hochster:1977}).  By taking a subsequence we may assume that the
sequence is non-increasing.   Let $u_t$ be an element of $R$ representing
the socle modulo $I_t$.  Then the injective hull of the residue field is $E = E_R(k) = \limt R/I_t$ and the image of
$u_t$ in $E$ is the socle element $u$ of $E$ for all $t$.  Moreover, because
the sequence is non-decreasing we may assume that for all $t$ there is
an injection $R/I_t \hookrightarrow R/I_{t+1}$  sending $u_t + I_t \mapsto u_{t+1} + I_{t+1}$.

Recall that a ring $R$ of positive prime characteristic is called \deff{$F$-finite} if the Frobenius
endomorphism  is module-finite.
Such rings are excellent \cite{Kunz:1976}, so if in addition $R$
is reduced then it is approximately Gorenstein.  Whenever  $R$ is reduced there 
is a well-defined ring of $q$th roots of $R$, denoted $\Rq$, which is
a finitely generated $R$-module for some (equivalently, all) $q$ precisely
when $R$ is $F$-finite.  In this case we will write $\Rq \cong R^{a_q} \oplus_R M_q$, where $M_q$ is a module with no free $R$ summands.  

The  characterization of the injective hull given above is very helpful in
proving the next Lemma, which shows how to compute the values of $a_q$
in a special case.  By $I\brq$ we mean the ideal $(i^q: i \in I)$.

\begin{lemma}\label{aq}
Let $(R,\m,k)$ be a reduced, $F$-finite ring with perfect residue field
$k$.  Then $a_q = \length_R(R/(I_t\brq : u_t^q))$ for all $t \gg 0$.
\end{lemma}

\begin{proof}
This result is a special case of Corollary 2.8 of \cite{Aberbach-Enescu:2003}.
However, we give a proof here for the benefit of the reader.  We will use the
 fact that over an approximately Gorenstein ring, a homomorphism  $f:R \to M$, where $M$ is finitely generated, has a splitting over $R$ if and only
if for all $t$, $f(u_t) \notin I_t M$ (see \cite{Hochster:1977}).  

Fix $q$, and write   $\Rq \cong R^{a_q} \oplus_R M_q$ as above.  We first claim that for $t \gg 0$, $u_t M_q \inc I_t M_q$, since for any minimal generator
of $M_q$, the map  $R x \to M_q$ does not split, and hence, $x u_t \in I_t M_q$.  The claim follows since $M_q$ is a finitely generated $R$-module.
We will also use the fact that if $I$ is an $\m$-primary ideal then
$\length_R(R/I\brq) = \length_R (\Rq/I \Rq)$, since $k$ is perfect.

Thus, for any $t \gg 0$, we have
\begin{align*}
\length(R/(I_t\brq : u_t^q))  &= \length(R/I_t\brq) - \length(R/(I_t,u_t)\brq)
= \length(\Rq/I_t \Rq) - \length(\Rq/ (I_t, u_t)\Rq)  \\
 &= \length(R^{a_q}/I_t R^{a_q}) + \length( M_q/I_t M_q) 
   -\left(\length(R^{a_q}/(I_t, u_t)R^{a_q}) + \length (M_q/(I_t,u_t)M_q)\right) \\
  &= a_q \cdot 1 + \length( M_q/I_t M_q) -  \length (M_q/(I_t,u_t)M_q) = a_q,
\end{align*}
since $(I_t, u_t)M_q = I_t M_q$ (for $t \gg 0$).
\end{proof}

We will need to pass to a $\Ga$ construction as described in \cite{HH:1994},
Section 6.  We refer the reader to \cite{HH:1994} for details.  What
we need to know is as follows.  Let $(R,\m, k)$ be a complete ring of characteristic
$p$.  Then $R \to R^\Ga$ is a faithfully flat, purely inseparable extension, the maximal ideal of
$R^\Ga$ is $\m R^\Ga$, and $R^\Ga$ is $F$-finite.  Note that if $I \inc R$ is
an irreducible $\m$-primary ideal of $R$ then $I R^\Ga$ is is also an
irreducible $\m R^\Ga$-primary ideal of $R^\Ga$.  Moreover, if $E_R(R/\m) = \limt R/I_t$, then $E_{R^\Ga}(R^\Ga/\m R^\Ga) = E_R(R/\m) \otimes_R R^\Ga = \limt R^\Ga/I_t R^\Ga$.

Our main theorem is 
\begin{thm}
Let $(R,\m,k)$ be an excellent local domain of positive prime characterstic.
Suppose that $\Tor_1(R^+, k) =0$.  Then $R$ is regular.
\end{thm}

\begin{proof}

By \cite{Schoutens:2003}, Theorem 1.2, the ring $R$ is weakly $F$-regular, therefore
a Cohen-Macaulay,  normal domain.  In particular, $R$ is approximately Gorenstein.   Also $R \to R^+$ is cyclically pure.
The assumption that $\Tor_1(R^+, k) =0$ and an induction on length shows
that for any $\m$-primary ideal $I \inc R$ and element $x$ we have $I R^+:_{R^+} x = (I:_R x)R^+$.

We first claim that for all $q$ and all $t$, $I_t\brq :_R u_t^q \subseteq \m\brq$.  To see this suppose that $v u_t^q \in I_t\brq$.  Taking
$q$th roots shows that $v^{1/q} \in I_t R^+:_{R^+} u_t = \m R^+$, and
hence that $v \in (\m\brq)^+ = \m\brq$ (by cyclic purity of $R$ in $R^+$).  This shows that 
for all $q$ and for all $t$, $\length(R/(I_t\brq:u_t^q)) \ge \length(R/\m\brq)$,
which is greater than or equal to $q^d$ (\cite{Kunz:1969}).

We consider $R \to \hat R \to (\hat R)^\Gamma = S$ for any
Gamma extension of $\hat R$.  In particular we may take $\Gamma$ to be the empty set, in which case
the residue field of $S$ is perfect.  Then by faithful flatness and the
fact that the maximal ideal of $S$ is $\m S$,
$\length_R(R/(I_t\brq:_R u_t^q)) = \length_S(S/(I_t S\brq:_S u_t^q))$.  
Since $u_t^q \notin I_t S\brq$ for all $t$, the ring $S$ is $F$-pure,
and hence reduced.
Thus by Lemma~\ref{aq},
for large enough $t$ (depending on $q$), $\length_R(R/(I_t\brq:_R u_t^q)) = a_q(S)$  is 
the number of $S$-free summands in $S^{1/q}$.    Since
$S$ has perfect residue field, the rank of $S^{1/q}$ as an $S$-module
is precisely $q^d$, hence $a_q(S) \le q^d$.  We have now shown that
$q^d \ge \length(R/(I_t\brq:u_t^q)) \ge \length(R/\m\brq) \ge q^d$
Thus $\length(R/\m\brq) = q^d$ and $R$ is regular \cite{Kunz:1969}.
\end{proof}

\end{document}